\documentclass[report]{llncs}
\usepackage{amsmath}
\usepackage{greekctr}
\usepackage{chngcntr}
\usepackage{booktabs} 
\usepackage{caption} 
\usepackage{subcaption} 
\usepackage{graphicx}
\usepackage{pgfplots}
\usepackage[all]{nowidow}
\usepackage[utf8]{inputenc}
\usepackage{tikz}
\usetikzlibrary{er,positioning,bayesnet}
\usepackage{multicol}
\usepackage{algpseudocode,algorithm,algorithmicx}
\usepackage{listings}
\usepackage{hyperref}
\usepackage[inline]{enumitem} 
\usepackage{theoremref}
\usepackage{amsfonts}
\usepackage{amssymb}
\usepackage{hyperref}

\usepackage{plain}
\usepackage{cleveref}

\definecolor{blue}{HTML}{1F77B4}
\definecolor{orange}{HTML}{FF7F0E}
\definecolor{green}{HTML}{2CA02C}
\pgfplotsset{compat=1.14}

\usepackage{hyperref}

\usepackage{cleveref}

\begin{document}
\title{Some Extensions of Endo-Noetherian Rings}
%
%
\author{R. M. Salem\inst{1} \and R. E. Abdel-Khalek\inst{1} \and N. Abdelnasser\inst{2}}

\institute{
Department of Mathematics, Faculty of Science, Al-Azhar University, Nasr City 11884, Cairo, Egypt\\
\email{rsalem\_02@hotmail.com, refaat\_salem@cic-cairo.com, \\ramy\_ama@yahoo.com, ramyabdel\_khalek@azhar.edu.eg}
\and
Mathematics Department, Faculty of Science, Assiut University, Assiut, Egypt\\
\email{neamanasser7@gmail.com}
}
%
%
%
\maketitle         
  \begin{abstract}


In this article, we proceed on the transfer of the left endo-Noetherian property on certain ring extensions. We transfer of the right (left) endo-Noetherian property to the right (left) quotient rings. For a subring $T$ of $R$ and a finite set of indeterminates $X$, we prove that $T + XR[[X]]$ is left endo-Noetherian if and only if $R[[X]]$ is left endo-Noetherian. In addition, we prove that the subring $\Lambda :=\{ f \in  R[[S,\omega ]]: f(1) \in T \}$ of the skew generalized power series ring $R[[S, \omega]]$ is left endo-Noetherian if and only if $R[[S, \omega]]$ is left endo-Noetherian. Also, we study the left endo-Noetherian property over the amalgamated duplication rings $R \bowtie I$ and $ R \bowtie ^f J$. Finally, we introduce additional results on left endo-Noetherian rings.

\end{abstract}

\section{Introduction}
Throughout this paper all rings are associative with identity (not necessarily commutative).
In 2009, A. Kaidi and E. Sanchez introduced the class of endo-Noetherian modules \cite{kaidimodules}. A left module $_{R}M$ of a ring $R$ is called endo-Noetherian if it satisfies the ascending chain condition for endomorphic kernels. A ring $R$ is called left endo-Noetherian if $_{R}R$ is endo-Noetherian as a left module. Equivalently, $R$ is left endo-Noetherian if the ascending chain of left annihilators $\ell .\mathrm {ann} _{R}  (r_1 )\subseteq  \ell .\mathrm {ann}
 _{R}  (r_2 )\subseteq \dots $ stabilizes for each sequence $(r_i)_{i\in \mathbb{N}}$ (i.e. there exists a positive integer $n$ such that $\ell .\mathrm {ann} _{R}  (r_{k})=\ell .\mathrm {ann} _{R}  (r_{n})$ for each $k\geq n$).
Similarly, $R$ is right endo-Noetherian if the ascending chain of right annihilators $ r .\mathrm {ann} _{R}  (r_1 )\subseteq  r .\mathrm {ann}
 _{R}  (r_2 )\subseteq \dots $ stabilizes for each sequence $(r_i)_{i\in \mathbb{N}}$. The class of endo-Noetherian lies between the
class of iso-Noetherian and the class of strongly hopfian. A right $R$-module $M$ is
iso-Noetherian if for every ascending chain ${\displaystyle \mathrm M_{1}\subseteq \mathrm M_{2} \subseteq \dots}$ of right submodules of $M$, there exists
an index $n \geq 1$ such that $M_{n} \simeq M_{i}$ for every $i \geq n$. A ring $R$ is called right iso-Noetherian
if the right $R$-module $R$ is iso-Noetherian and $R$ is called right strongly Hopfian if for every $a \in R$ there exists a positive integer $n$ such that  ${\displaystyle  r .\mathrm {ann}(a^{n})=r .\mathrm {ann}(a^{n+1})}$). Also, every Noetherian rings is endo-Noetherian but the converse is not true. These relations and some counter examples are shown in \cite{mohamed2023endo}. In general, the submodules of endo-Noetherian modules need not be endo-Noetherian, see \cite{kaidimodules}.
 In \cite{gouaid2020endo}, Gouaid et al. studied the endo-Noetherian property with quotient rings in the commutative case. They gave an example of a commutative ring $R$ and a multiplicative subset $S$ of
$R$ such that the localization $R_{S}$ of $R$ is Noetherian (so endo-Noetherian) but $R$ is not endo-Noetherian. Also, they introduced
a sufficient condition for  $R_{S}$ satisfies the endo-Noetherian property
implies that $R$ is endo-Noetherian. In \cite{kaidimodules}, Kaidi gave an example to show that the quotients of endo-Noetherian modules need not be endo-Noetherian. 
\\The purpose of this paper is to study the left endo-Noetherian property on some ring extensions. In Section 2, we transfer of the right (left) endo-Noetherian property to the right (left) quotient rings.
In section 3, we prove that $T + XR[[X]]$ is left endo-Noetherian if and only if $R[[X]]$ is left endo-Noetherian, for a subring $T$ of $R$ and a finite set of indeterminates $X$. In addition,  We introduce the structure $\Lambda :=\{ f \in  R[[S,\omega ]]: f(1) \in T \}$ which is a subring of the skew generalized power series $ R[[S,\omega ]]$. We prove that the subring $\Lambda :=\{ f \in  R[[S,\omega ]]: f(1) \in T \}$ of the skew generalized power series ring $R[[S, \omega]]$ is left endo-Noetherian if and only if $R[[S, \omega]]$ is left endo-Noetherian. 

Let us recall the following
notion. Let $S=(R_{n})_{n\in \mathbb{N}}$ be an increasing sequence of rings, $ R=\cup _{n\in \mathbb{N}} R_{n} $, their union and let $S[x]$ be the ring of polynomials with coefficients of degree $n$ in $R_{n}$. Then $S[x]$ is a subring of the ring of polynomials
$R[x]$, see \cite{ahmed2015s}. In \cite[Corollary 3]{mohamed2023endo}, the authors introduced the equivalent conditions for the polynomial rings over an Armendariz ring to be left endo-Noetherian. In Section 4, we introduce the equivalent conditions for the structure $S[x]$ to be left endo-Noetherian. Finally, we study when the amalgamated
duplication $ R \bowtie I $ and $ R \bowtie ^f J$ satisfy the left endo-Noetherian property.

\section{Transfer of the Endo-Noetherian Property to the Quotient Rings}
\begin{definition}\cite{kaidimodules}
   A ring $R$ is called left endo-Noetherian if the ascending chain of left annihilators $\ell .\mathrm {ann} _{R}  (r_1 )\subseteq  \ell .\mathrm {ann}
 _{R}  (r_2 )\subseteq \dots $ stabilizes for each sequence $(r_i)_{i\in \mathbb{N}}$ (i.e. there exists a positive integer $n$ such that $\ell .\mathrm {ann} _{R}  (r_{k})=\ell .\mathrm {ann} _{R}  (r_{n})$ for each $k\geq n$).
\end{definition}
\begin{definition}\cite[2.1.13]{mcconnell2001noncommutative} A multiplicatively closed subset
$S$ of a ring $R$ 
is said to be a left
Ore set if for each $ r \in R$ and $s \in S$ there exists $r^{'} \in R$, $s^{'} \in S$ such that $rs^{'}  = sr^{'}$ (i.e. $Sr \cap  Rs \neq \phi $).
\end{definition}
  
 Unlike commutative rings, the existence of a right (or left) quotient ring is not assured for noncommutative rings. Furthermore, one-sidedness (right or left) does not necessarily indicate the presence of the other (see \cite[p. 45]{mcconnell2001noncommutative}). We denote the left quotient ring by $Q$ and the right quotient ring by $Q^{'}$.
In this section, we examine how the right endo-Noetherian property is transferred from the ground ring $R$ to the right quotient ring $Q^{'}$ and vice versa.

\begin{proposition}
    Let $R$ be a ring and $S$ a right Ore set consists of regular elements. If the right quotient ring $Q^{'}$ is right endo-Noetherian, then $R$ is also right endo-Noetherian.
\end{proposition}
\begin{proof}
Assume that $Q^{'}$ is right endo-Noetherian and $(r_k )_{(k \in N)}$ is a sequence of elements of $R$ such that $I_1\subseteq I_2\subseteq \dots $ in $R$, where $ I_i= r .\mathrm {ann} _{R}(r_i)$ is a right ideal in $R$. By \cite[Proposition 1.16]{mcconnell2001noncommutative}, $I_i Q^{'}=\{ x s^{-1}  \mid x\in I_i, s\in S \} $ is a right ideal in $Q^{'}$ for each  $i\in \mathbb{N}$. \\
 One can easily check that $IQ^{'}= r .\mathrm {ann} _{Q^{'}}(r)$, where  $I= r .\mathrm {ann} _{R}(r)$. Let $x \in IQ^{'}$. Then there exist $i\in I$, $s_1 \in S$ such that $x=i{s_1}^{-1}$, and $ri=0$. Thus $rx=r(i{s_1}^{-1})=(ri){s_1}^{-1}=0$, and $x\in r .\mathrm {ann} _{Q^{'}}(r)$. Also let $ r^{'} s^{'-1 }\in r .\mathrm {ann} _{Q^{'}}(r)$. Then $rr^{'}=rr^{'} s^{'-1 }=0$, and $r^{'} \in r .\mathrm {ann} _{R}(r) =I$. Hence $r^{'} s^{'-1 }\in IQ^{'}$. \\
 We will show that $I_j Q^{'}\subseteq I_{j+1} Q^{'} $ for each $j \in \mathbb{N}$. Let $x\in I_j Q^{'}$. Then there exist $i\in I_j$, $s\in S$ such that $x=is^{-1}$. Since $i\in I_j \subseteq I_{j+1}=r .\mathrm {ann} _{R}(r_{j+1})$, $ r_{j+1} i=0$. Where $s^{-1} \in Q^{'} $, we have $r_{j+1} is^{-1}=0$, and $is^{-1}= x\in r .\mathrm {ann}_{Q^{'}} (r_{j+1})=I_{j+1}Q^{'}$. Now, since $Q^{'}$ is right endo-Noetherian, there exists a positive integer $n$ such that $I_{k} Q^{'}=I_{n} Q^{'}$ for each $k\geq n$. We will show that $I_{k}=I_{n}$. Let $r \in I_{k}=r .\mathrm {ann} _{R}(r_{k})$. Then $ r_{k} r s^{-1}=r_{k} r=0$. Since $s^{-1}\in Q^{'}$, we have $r s^{-1} \in r .\mathrm {ann}_{Q^{'}} (r_k )=I_k Q^{'}=I_n Q^{'}$. Therefore $r_{n} r= r_{n} r s^{-1}=0$, and $ r \in r .\mathrm {ann} _{R}(r_n )=I_n$. Hence $R$ is right endo-Noetherian.
\end{proof}
\begin{proposition}
   Let $R$ be a ring and $S$ a right Ore set consists of regular elements. If $R$ is left endo-Noetherian, then $Q^{'}$ is also left endo-Noetherian.
\end{proposition}
\begin{proof}
    Assume that $R$ is left endo-Noetherian and $(r_i {s_{i}}^{-1} )_{i\in \mathbb{N}}$ is a sequence of elements of $Q^{'}$ such that $B_1\subseteq B_2\subseteq \dots $ in $Q^{'}$ where $B_i=\ell .\mathrm {ann} _{Q^{'}}(r_i {s_{i}}^{-1})$ is a left ideal of $Q^{'}$. By \cite[Proposition 1.16]{mcconnell2001noncommutative}, $B_{i}\cap R$ is a left  ideal of $R$, where $B_{i}\cap R=\{ a_i\in R \mid a_i 1^{-1}\in B_i\} $  for each  $ i\in \mathbb{N}$.\\
    One can easily check that  $B_{i}\cap R=\ell .\mathrm {ann} _{R}(r_i)$. Let $b\in B_{i}\cap R$. Then $b1^{-1}\in B_i =\ell .\mathrm {ann} _{Q^{'}}(r_i {s_{i}}^{-1})$, and $b1^{-1} r_i {s_{i}}^{-1}=0$. Since ${s_{i}}^{-1}$ is a unit in $Q^{'}$, we have $br_i=0$, and $b\in \ell .\mathrm {ann} _{R}(r_i)$. Also, let $b\in \ell .\mathrm {ann} _{R}(r_i)$. Thus $br_i {s_{i}}^{-1}=br_i=0$. Therefore $b1^{-1}\in \ell .\mathrm {ann} _{Q^{'}}(r_{i} {s_{i}}^{-1})=B_i$.\\
     We will show that $B_i\cap R \subseteq B_{i+1}\cap R$ for each  $i\in \mathbb{N}$. Let $x_i\in B_i\cap R$. Then $x_i 1^{-1}\in B_i\subseteq B_{i+1}$, and $x_i\in B_{i+1}\cap R$.\\
     Now, since $R$ is left endo-Noetherian, there exists a positive integer $n$ such that $B_k\cap R=B_n\cap R$ for each $k\geq n$. By \cite[Proposition 1.16]{mcconnell2001noncommutative}, $B_k=(B_k\cap R)Q^{'}$ is the ideal which generated by $B_k\cap R$. Also $B_n=(B_n\cap R)Q^{'}$ is the ideal which generated by $B_n \cap R$. Therefore $B_k=B_n$.  Hence $Q^{'}$ is left endo-Noetherian.
     \end{proof}

\begin{remark}\label{rem1}
    Let $R$ be a ring and $S\subseteq R$ an Ore set consists of regular elements. Then from \cite[Theorem 2.1.12]{mcconnell2001noncommutative}, $R$ has a left quotient ring $Q$ together with a ring homomorphism $f: R\longrightarrow Q$ and a right quotient ring $Q^{'}$  together with a ring homomorphism $f^{'}: R\longrightarrow Q^{'}$. Also from \cite[Corollary 2.1.4]{mcconnell2001noncommutative}, we have $Q\cong Q^{'}$. It is possible to find that the ring isomorphism $\phi: Q\longrightarrow Q^{'}$ defined as follows
 $\phi(f(s)^{-1} f(r))=f^{'} (s)^{-1} f^{'} (r)$, where $r\in R$, $s\in S$.

\end{remark}
   In the following theorem, we use another way to prove that the ground ring $R$ is left endo-Noetherian if and only if the right quotient ring $Q^{'}$ is.

\begin{theorem}
   Let $R$ be a ring and $S$ an Ore set consists of regular elements. Then the following assertions are equivalent:
   \begin{enumerate}
       \item $R$ is right endo-Noetherian.
       \item $Q$ is right endo-Noetherian.
   \end{enumerate}
\end{theorem}
\begin{proof}
$(a) \Longrightarrow (b)$. Let $\left( f(s_{i})^{-1} f(r_i) \right)_{i \in \mathbb{N}}$ 
be a sequence of elements of $Q$ for some $r_i \in R$, $s_i \in S$ such that:
\[
r \cdot \mathrm{ann}_Q(f(s_{1})^{-1} f(r_1)) 
\subseteq r \cdot \mathrm{ann}_Q(f(s_{2})^{-1} f(r_2)) 
\subseteq \cdots.
\]

We will show that:
\[
r \cdot \mathrm{ann}_R(r_i) \subseteq r \cdot \mathrm{ann}_R(r_{i+1}), \quad \text{for each } i \in \mathbb{N}.
\]

Let $b \in r \cdot \mathrm{ann}_R(r_i)$, i.e., $r_i b = 0$. Then:
\[
f(r_i b) = 0 \Rightarrow f(r_i) f(b) = 0 \Rightarrow f(s_i)^{-1} f(r_i) f(b) = 0.
\]

Hence:
\[
f(b) \in r \cdot \mathrm{ann}_Q(f(s_i)^{-1} f(r_i)) 
\subseteq r \cdot \mathrm{ann}_Q(f(s_{i+1})^{-1} f(r_{i+1})).
\]

Thus:
\[
f(s_{i+1})^{-1} f(r_{i+1}) f(b) = 0 
\Rightarrow f(r_{i+1}) f(b) = 0 
\Rightarrow f(r_{i+1} b) = 0.
\]

So $r_{i+1} b \in \ker f$. Since $S$ consists of regular elements, and $\mathrm{ass}(S) = 0 = \ker f$, it follows that $r_{i+1} b = 0$. Therefore:
\[
b \in r \cdot \mathrm{ann}_R(r_{i+1}).
\]

Now, since $R$ is right endo-Noetherian, there exists a positive integer $n$ such that:
\[
r \cdot \mathrm{ann}_R(r_k) = r \cdot \mathrm{ann}_R(r_n) \quad \text{for all } k \geq n.
\]

Let $f(s)^{-1} f(r) \in r \cdot \mathrm{ann}_Q(f(s_k)^{-1} f(r_k))$, so:
\[
f(s_k)^{-1} f(r_k) f(s)^{-1} f(r) = 0.
\]

Since $S$ is an Ore set consisting of regular elements, it follows from Remark 1 that $R$ has a right quotient ring $Q' \cong Q$, with an isomorphism:
\[
\phi: Q \longrightarrow Q'
\]
such that:
\[
\phi(f(s)^{-1} f(r)) = f'(s)^{-1} f'(r).
\]

\begin{align*}
\phi(f(s_k)^{-1} f(r_k) f(s)^{-1} f(r)) &= 0, \\
\phi(f(s_k)^{-1} f(r_k)) \cdot \phi(f(s)^{-1} f(r)) &= 0, \\
(f'(s_k)^{-1} f'(r_k)) (f'(s)^{-1} f'(r)) &= 0.
\end{align*}

Since $f'(s_k)^{-1}$ is a unit in $Q'$, we have
\[
f'(r_k)(f'(s)^{-1} f'(r)) = 0.
\]
Since $f'(s)^{-1} f'(r) \in Q'$, we can write
\[
f'(s)^{-1} f'(r) = f'(r') f'(s')^{-1}
\]
with $r' \in R$, $s' \in S$. Then
\[
f'(r_k) (f'(r') f'(s')^{-1}) = 0.
\]
Since $f'(s')^{-1}$ is a unit in $Q'$, ...

We have:
\[
f'(r_k) f'(r') = 0 \Rightarrow f'(r_k r') = 0 \Rightarrow r_k r' \in \ker f' = \{0\} = \mathrm{ass}(S) \Rightarrow r_k r' = 0.
\]

Hence:
\[
r' \in r \cdot \mathrm{ann}_R(r_k) = r \cdot \mathrm{ann}_R(r_n),
\]
so:
\[
r_n r' = 0 \Rightarrow f'(r_n r') = 0 \Rightarrow f'(r_n) f'(r') = 0.
\]

Then:
\[
f'(s_n)^{-1} f'(r_n) f'(r') f'(s')^{-1} = 0.
\]

But since:
\[
f'(r') f'(s')^{-1} = f'(s)^{-1} f'(r),
\]
we get:
\[
f'(s_n)^{-1} f'(r_n) f'(s)^{-1} f'(r) = 0.
\]

Therefore:
\[
\phi(f(s_n)^{-1} f(r_n)) \cdot \phi(f(s)^{-1} f(r)) = 0.
\]

Since $\phi$ is an isomorphism, it follows that:
\[
f(s_n)^{-1} f(r_n) f(s)^{-1} f(r) = 0,
\]
so:
\[
f(s)^{-1} f(r) \in r \cdot \mathrm{ann}_Q(f(s_n)^{-1} f(r_n)).
\]

Hence, $Q$ is right endo-Noetherian.

\noindent
\textbf{$(b)\Longrightarrow (a)$}. Assume that $Q$ is right endo-Noetherian, and let $(r_k)_{k \in \mathbb{N}}$ be a sequence in $R$ such that:
\[
r \cdot \mathrm{ann}_R(r_1) \subseteq r \cdot \mathrm{ann}_R(r_2) \subseteq \cdots
\]

We will show that:
\[
r \cdot \mathrm{ann}_Q(f(s_0)^{-1}f(r_i)) \subseteq r \cdot \mathrm{ann}_Q(f(s_0)^{-1}f(r_{i+1}))
\]
for some $s_0 \in S$ and for each $i \in \mathbb{N}$.

Let $f(s)^{-1}f(r) \in r \cdot \mathrm{ann}_Q(f(s_0)^{-1}f(r_i))$, so:
\[
f(s_0)^{-1} f(r_i) f(s)^{-1} f(r) = 0.
\]

As above, the left quotient ring $Q$ is isomorphic to the right quotient ring $Q'$ via an isomorphism:
\[
\phi : Q \longrightarrow Q' \quad \text{such that} \quad \phi(f(s)^{-1}f(r)) = f'(s)^{-1}f'(r).
\]

Then:
\begin{align*}
0 &= \phi(f(s_0)^{-1}f(r_i)f(s)^{-1}f(r)) \\
  &= \phi(f(s_0)^{-1}f(r_i)) \cdot \phi(f(s)^{-1}f(r)) \\
  &= f'(s_0)^{-1}f'(r_i) \cdot f'(s)^{-1}f'(r).
\end{align*}

Since $f'(s_0)^{-1}$ is a unit in $Q'$, we get:
\[
f'(r_i) \cdot f'(s)^{-1}f'(r) = 0.
\]

Now, since $f'(s)^{-1}f'(r) \in Q'$, we can write:
\[
f'(s)^{-1}f'(r) = f'(r') f'(s')^{-1}, \quad \text{for some } r' \in R, \; s' \in S.
\]

Then:
\[
f'(r_i)f'(r')f'(s')^{-1} = 0 \Rightarrow f'(r_i r') = 0,
\]
and hence $r_i r' \in \ker f'$.

Since $S$ consists of regular elements, we have $\ker f' = \operatorname{ass} S = 0$, so:
\[
r_i r' = 0 \Rightarrow r' \in r \cdot \mathrm{ann}_R(r_i) \subseteq r \cdot \mathrm{ann}_R(r_{i+1}) \Rightarrow r_{i+1} r' = 0.
\]

Thus:
\[
f'(r_{i+1} r') = 0 \Rightarrow f'(r_{i+1}) f'(r') = 0.
\]

Now, multiplying both sides:
\[
f'(s_0)^{-1} f'(r_{i+1}) f'(r') f'(s')^{-1} = 0.
\]

But since $f'(r') f'(s')^{-1} = f'(s)^{-1} f'(r)$, we have:
\[
\phi(f(s_0)^{-1} f(r_{i+1})) \cdot \phi(f(s)^{-1} f(r)) = 0.
\]

Using that $\phi$ is an isomorphism, it follows that:
\[
f(s_0)^{-1} f(r_{i+1}) f(s)^{-1} f(r) = 0,
\]
so:
\[
f(s)^{-1} f(r) \in r \cdot \mathrm{ann}_Q(f(s_0)^{-1} f(r_{i+1})).
\]

Therefore:
\[
r \cdot \mathrm{ann}_Q(f(s_0)^{-1} f(r_i)) \subseteq r \cdot \mathrm{ann}_Q(f(s_0)^{-1} f(r_{i+1})).
\]

Now, since $Q$ is right endo-Noetherian, there exists $n \in \mathbb{N}$ such that:
\[
r \cdot \mathrm{ann}_Q(f(s_0)^{-1} f(r_k)) = r \cdot \mathrm{ann}_Q(f(s_0)^{-1} f(r_n)) \quad \text{for all } k \geq n.
\]

Let $\alpha \in r \cdot \mathrm{ann}_R(r_k)$, i.e., $r_k \alpha = 0 \Rightarrow f(r_k \alpha) = 0$.

Hence:
\[
f(r_k) f(\alpha) = 0 \Rightarrow f(s_0)^{-1} f(r_k) f(\alpha) = 0.
\]

So:
\[
f(\alpha) \in r \cdot \mathrm{ann}_Q(f(s_0)^{-1} f(r_k)) = r \cdot \mathrm{ann}_Q(f(s_0)^{-1} f(r_n)) \Rightarrow f(s_0)^{-1} f(r_n) f(\alpha) = 0.
\]

Since $f(s_0)^{-1}$ is a unit:
\[
f(r_n) f(\alpha) = 0 \Rightarrow f(r_n \alpha) = 0 \Rightarrow r_n \alpha \in \ker f = 0.
\]

Thus $r_n \alpha = 0$, and so $\alpha \in r \cdot \mathrm{ann}_R(r_n)$, hence $R$ is right endo-Noetherian.
           \end{proof}
    
\section{Endo-Noetherian Rings of The 
Form $T+XR[[X]]$ and Its Related Rings}

 In this section, we examine the endo-Noetherian property on a particular subring  of the formal power series ring $R[[X]]$, such as the subring $T + XR[[X]]$, where $X:=\{x_1,x_2,...,x_n\}$ is a finite set of indeterminate and $T$ is a subring of $R$. However, we generalize \cite[Proposition 2.1]{hamed2021rings} in the following theorem.
\begin{theorem}
 Let $T\subseteq R$ be an extension of rings. Then the following conditions are equivalent:
 \begin{enumerate}
     \item $T+XR[[X]]$ is left endo-Noetherian. 
      \item $R[[X]]$ is left endo-Noetherian.
 \end{enumerate}
\end{theorem}
\begin{proof}
\textbf{($a \Rightarrow b$).} Let $(f_i)_{i \in \mathbb{N}}$ be a sequence in $R[[X]]$,
\begin{equation*}
    f_i = \sum_{i_1, i_2, \dots, i_n = 0}^{\infty} b_{i_1, i_2, \dots, i_n} x_1^{i_1} x_2^{i_2} \dots x_n^{i_n},
\end{equation*}
such that
\[
\ell.\mathrm{ann}_{R[[X]]}(f_1) \subseteq \ell.\mathrm{ann}_{R[[X]]}(f_2) \subseteq \cdots.
\]

Since $x_1 f_i \in T + X R[[X]]$, we show:
\[
\ell.\mathrm{ann}_{T + X R[[X]]}(x_1 f_1) \subseteq \ell.\mathrm{ann}_{T + X R[[X]]}(x_1 f_2) \subseteq \cdots.
\]

Let $q \in \ell.\mathrm{ann}_{T + X R[[X]]}(x_1 f_1)$. Then:
\[
q f_1 = q x_1 f_1 = 0,
\]
and thus $q \in \ell.\mathrm{ann}_{R[[X]]}(f_1) \subseteq \ell.\mathrm{ann}_{R[[X]]}(f_2)$, so:
\[
q f_2 = q x_1 f_2 = 0 \Rightarrow q \in \ell.\mathrm{ann}_{T + X R[[X]]}(x_1 f_2).
\]

Now, since $T + X R[[X]]$ is left endo-Noetherian, there exists $n \in \mathbb{N}$ such that:
\[
\ell.\mathrm{ann}_{T + X R[[X]]}(x_1 f_k) = \ell.\mathrm{ann}_{T + X R[[X]]}(x_1 f_n) \quad \text{for all } k \geq n.
\]

We show:
\[
\ell.\mathrm{ann}_{R[[X]]}(f_k) = \ell.\mathrm{ann}_{R[[X]]}(f_n) \quad \text{for all } k \geq n.
\]

Let $g \in \ell.\mathrm{ann}_{R[[X]]}(f_k)$. Then:
\[
g f_k = x_1 g x_1 f_k = 0 \Rightarrow x_1 g \in \ell.\mathrm{ann}_{T + X R[[X]]}(x_1 f_k)
\subseteq \ell.\mathrm{ann}_{T + X R[[X]]}(x_1 f_n).
\]

Thus:
\[
g f_n = x_1 g x_1 f_n = 0 \Rightarrow g \in \ell.\mathrm{ann}_{R[[X]]}(f_n).
\]

Hence, $R[[X]]$ is left endo-Noetherian.

\vspace{1em}
\noindent
\textbf{($b \Rightarrow a$).} Let $(q_i)_{i \in \mathbb{N}}$ be a sequence in $T + X R[[X]]$ such that:
\[
\ell.\mathrm{ann}_{T + X R[[X]]}(q_1) \subseteq \ell.\mathrm{ann}_{T + X R[[X]]}(q_2) \subseteq \cdots.
\]

We show:
\[
\ell.\mathrm{ann}_{R[[X]]}(q_i) \subseteq \ell.\mathrm{ann}_{R[[X]]}(q_{i+1}) \quad \text{for each } i \in \mathbb{N}.
\]

Let $g \in \ell.\mathrm{ann}_{R[[X]]}(q_i)$. Then:
\[
g q_i = x_1 g q_i = 0 \Rightarrow x_1 g \in \ell.\mathrm{ann}_{T + X R[[X]]}(q_i)
\subseteq \ell.\mathrm{ann}_{T + X R[[X]]}(q_{i+1}).
\]

Thus:
\[
g q_{i+1} = x_1 g q_{i+1} = 0 \Rightarrow g \in \ell.\mathrm{ann}_{R[[X]]}(q_{i+1}).
\]

Since $R[[X]]$ is left endo-Noetherian, there exists $n \in \mathbb{N}$ such that:
\[
\ell.\mathrm{ann}_{R[[X]]}(q_k) = \ell.\mathrm{ann}_{R[[X]]}(q_n) \quad \text{for all } k \geq n.
\]

We now show:
\[
\ell.\mathrm{ann}_{T + X R[[X]]}(q_k) = \ell.\mathrm{ann}_{T + X R[[X]]}(q_n) \quad \text{for all } k \geq n.
\]

Let $q \in \ell.\mathrm{ann}_{T + X R[[X]]}(q_k)$. Then $q q_k = 0$, and since:
\[
q \in \ell.\mathrm{ann}_{R[[X]]}(q_k) \subseteq \ell.\mathrm{ann}_{R[[X]]}(q_n) \Rightarrow q q_n = 0,
\]
we conclude:
\[
q \in \ell.\mathrm{ann}_{T + X R[[X]]}(q_n).
\]

Hence, $T + X R[[X]]$ is left endo-Noetherian.
\end{proof}

To show that the subring $ \Lambda$ that corresponds $T+xR[[x]]$ of the form $\{ f \in  R[[S,\omega ]]: f(1) \in T \}$ is left endo Noetherian if and only if $ R[[S,\omega ]]$ is left endo-Noetherian the following proposition is essential:

\begin{proposition}\cite[Proposition 4.2.]{dina2024almost}\label{prop22}
    Let $R$ be a ring, $(S, \preceq)$ a totally ordered monoid, $\omega:S\longrightarrow End(R)$ a monoid homomorphism, and $R$ is $S$-compatible. Assume that for every $f \in R[[S,\omega ]] $, there exists $s_{0} \in supp{f}$. If $f(s_{0})$ is right (left) regular, then $f$ is right (left) regular.
\end{proposition}
From this proposition, we can determine a regular element in $ R[[S,\omega ]]$ as follows.
\begin{lemma}\label{prop3}
    Let $R$ be a ring, $(S,\preceq)$ a strictly ordered monoid satisfying the condition that $s \geq 1$ for every $s\in S$, and $\omega:S\longrightarrow End(R)$ a monoid homomorphism. Assume that $R$ is $S$-compatible. Then $e_{s}$ is a regular element in $ R[[S,\omega ]]$.
\end{lemma}
Now, we can conclude the main result of this section as follows.
\begin{theorem}\label{them11}
    Let $T \subseteq R$ be an extension of rings, $(S, \preceq)$ a strictly ordered monoid satisfying the condition that $s \geq 1$ for every $s\in S$,  $\omega:S\longrightarrow End(R)$ a monoid homomorphism and $\Lambda :=\{ f \in  R[[S,\omega ]]: f(1) \in T \}$ a subring of $ R[[S,\omega ]]$. Assume that $R$ is $S$-compatible, then the following conditions are equivalent:
     \begin{enumerate}
         \item $ \Lambda$ is left endo-Noetherian.
         \item $  R[[S,\omega ]]$ is left endo-Noetherian.
     \end{enumerate}
\end{theorem}
\begin{proof}
\noindent
$(a)\Longrightarrow (b)$. Let $(f_i)_{i \in \mathbb{N}}$ be a sequence of elements of $R[[S, \omega]]$ such that
\[
\ell.\mathrm{ann}_{R[[S, \omega]]}(f_1)
\subseteq \ell.\mathrm{ann}_{R[[S, \omega]]}(f_2)
\subseteq \cdots.
\]

Since for $1 \ne s \in S$, we have $f_i e_s \in \Lambda$ for each $i \in \mathbb{N}$, we will show that
\[
\ell.\mathrm{ann}_\Lambda(f_i e_s) \subseteq \ell.\mathrm{ann}_\Lambda(f_{i+1} e_s)
\quad \text{for each } i \in \mathbb{N}.
\]

Let $h \in \ell.\mathrm{ann}_\Lambda(f_i e_s)$. Then
\[
h f_i e_s = 0.
\]
By Lemma 1, $e_s$ is a regular element in $R[[S, \omega]]$, so
\[
h f_i = 0 \quad \Rightarrow \quad h \in \ell.\mathrm{ann}_{R[[S, \omega]]}(f_i)
\subseteq \ell.\mathrm{ann}_{R[[S, \omega]]}(f_{i+1}).
\]

Hence,
\[
h f_{i+1} e_s = h f_{i+1} = 0
\quad \Rightarrow \quad h \in \ell.\mathrm{ann}_\Lambda(f_{i+1} e_s).
\]

Now, since $\Lambda$ is left endo-Noetherian, there exists a positive integer $n$ such that for all $k \geq n$:
\[
\ell.\mathrm{ann}_\Lambda(f_k e_s) = \ell.\mathrm{ann}_\Lambda(f_n e_s).
\]

We will show that
\[
\ell.\mathrm{ann}_{R[[S, \omega]]}(f_k) = \ell.\mathrm{ann}_{R[[S, \omega]]}(f_n)
\quad \text{for each } k \geq n.
\]

Let $g \in \ell.\mathrm{ann}_{R[[S, \omega]]}(f_k)$. Then
\[
e_s g f_k e_s = 0.
\]
Since $e_s g \in \Lambda$, and
\[
e_s g \in \ell.\mathrm{ann}_\Lambda(f_k e_s)
\subseteq \ell.\mathrm{ann}_\Lambda(f_n e_s),
\]
we have:
\[
e_s g f_n e_s = 0.
\]
Since $e_s$ is a regular element in $R[[S, \omega]]$, it follows that
\[
g f_n = 0 \quad \Rightarrow \quad g \in \ell.\mathrm{ann}_{R[[S, \omega]]}(f_n).
\]

Hence, $R[[S, \omega]]$ is left endo-Noetherian.

\vspace{0.5em}

\noindent
$(b)\Longrightarrow (a)$. Let $(q_i)_{i \in \mathbb{N}}$ be a sequence of elements of $\Lambda$ such that
\[
\ell.\mathrm{ann}_\Lambda(q_1) \subseteq \ell.\mathrm{ann}_\Lambda(q_2) \subseteq \cdots.
\]

We will show that
\[
\ell.\mathrm{ann}_{R[[S, \omega]]}(q_i) 
\subseteq 
\ell.\mathrm{ann}_{R[[S, \omega]]}(q_{i+1})
\quad \text{for each } i \in \mathbb{N}.
\]

Let $g \in \ell.\mathrm{ann}_{R[[S,\omega]]}(q_i)$. Then
\[
e_s g q_i = 0.
\]
Since $e_s g \in \Lambda$, we have:
\[
e_s g \in \ell.\mathrm{ann}_{\Lambda}(q_i)
\subseteq \ell.\mathrm{ann}_{\Lambda}(q_{i+1}),
\]
which implies:
\[
e_s g q_{i+1} = 0.
\]

Since $e_s$ is a regular element in $R[[S,\omega]]$, it follows that:
\[
g q_{i+1} = 0 \quad \Rightarrow \quad g \in \ell.\mathrm{ann}_{R[[S,\omega]]}(q_{i+1}).
\]

Now, since $R[[S,\omega]]$ is left endo-Noetherian, there exists a positive integer $n$ such that for all $k \geq n$:
\[
\ell.\mathrm{ann}_{R[[S,\omega]]}(q_k) = \ell.\mathrm{ann}_{R[[S,\omega]]}(q_n).
\]

We will show that:
\[
\ell.\mathrm{ann}_{\Lambda}(q_k) = \ell.\mathrm{ann}_{\Lambda}(q_n) \quad \text{for each } k \geq n.
\]

Let $q \in \ell.\mathrm{ann}_{\Lambda}(q_k)$. Then:
\[
q q_k = 0,
\]
and since:
\[
q \in \ell.\mathrm{ann}_{R[[S,\omega]]}(q_k)
\subseteq \ell.\mathrm{ann}_{R[[S,\omega]]}(q_n),
\]
we get:
\[
q q_n = 0 \quad \Rightarrow \quad q \in \ell.\mathrm{ann}_{\Lambda}(q_n).
\]

Hence, $\Lambda$ is left endo-Noetherian.

\end{proof}
If we assume that $\omega$ is the identity endomorphism, we have the following corollary.
\begin{corollary}
     Let $T $, $R$, $S$ be as in Theorem 3 and  $\Lambda :=\{ f \in R[[S]]: f(1) \in T \}$ a subring of $R[[S]]$. Then 
     \begin{enumerate}
         \item $ \Lambda$ is left endo-Noetherian if and only if $ R[[S]]$ is left endo-Noetherian.
        
     \item $T+xR[[x]]$ is left endo-Noetherian if and only if $R[[x]]$ is left endo-Noetherian.

     \item $T+xR[x]$ is left endo-Noetherian if and only if $R[x]$ is left endo-Noetherian.
 \end{enumerate}
\end{corollary}
It is well known if $R$ is $\sigma$-compatible then $\sigma$ is an injective homomorphism. The purpose of the following two propositions is to prove when $T+R[x,\sigma]x$ and $T+xR[x,\sigma]$ are respectively left endo-Noetherian and right endo-Noetherian.
\begin{proposition}
    Let $T \subseteq R$ be an extension of rings and $\sigma$ an injective endomorphism of $R$. Then the following conditions are equivalent:
    \begin{enumerate}
        \item $T+R[x,\sigma]x$ is left endo-Noetherian.
        \item $R[x,\sigma]$ is left endo-Noetherian.
    \end{enumerate}
   \end{proposition}
	\begin{proof}
\noindent
$(a)\Longrightarrow (b)$. Let $(f_{k})_{k \in \mathbb{N}}$ be a sequence of elements of $R[x, \sigma]$ such that
\[
\ell.\mathrm{ann}_{R[x,\sigma]}(f_1) \subseteq \ell.\mathrm{ann}_{R[x,\sigma]}(f_2) \subseteq \cdots.
\]
Since $f_k x \in T + R[x, \sigma]x$ for each $f_k \in R[x, \sigma]$, we will show that:
\[
\ell.\mathrm{ann}_{T + R[x,\sigma]x}(f_i x) \subseteq \ell.\mathrm{ann}_{T + R[x,\sigma]x}(f_{i+1} x)
\quad \text{for each } i \in \mathbb{N}.
\]

Let $q \in \ell.\mathrm{ann}_{T + R[x,\sigma]x}(f_i x)$. Then:
\[
q f_i = q f_i x = 0,
\]
and so $q \in \ell.\mathrm{ann}_{R[x,\sigma]}(f_i) \subseteq \ell.\mathrm{ann}_{R[x,\sigma]}(f_{i+1})$. Hence,
\[
q f_{i+1} x = q f_{i+1} = 0,
\]
and thus $q \in \ell.\mathrm{ann}_{T + R[x,\sigma]x}(f_{i+1} x)$.

Now, since $T + R[x, \sigma]x$ is left endo-Noetherian, there exists a positive integer $n$ such that for each $k \geq n$:
\[
\ell.\mathrm{ann}_{T + R[x,\sigma]x}(f_k x) = \ell.\mathrm{ann}_{T + R[x,\sigma]x}(f_n x).
\]

We will now show that:
\[
\ell.\mathrm{ann}_{R[x,\sigma]}(f_k) = \ell.\mathrm{ann}_{R[x,\sigma]}(f_n)
\quad \text{for each } k \geq n.
\]

Let $g \in \ell.\mathrm{ann}_{R[x,\sigma]}(f_k)$. Then:
\[
x g f_k x = g f_k = 0,
\]
and so $x g \in \ell.\mathrm{ann}_{T + R[x,\sigma]x}(f_k x) \subseteq \ell.\mathrm{ann}_{T + R[x,\sigma]x}(f_n x)$.

Thus,
\[
x g f_n x = \sigma(g) \sigma(f_n) x^2 = 0 \Rightarrow \sigma(g) \sigma(f_n) = 0.
\]
Since $\sigma$ is injective, we conclude $g f_n = 0$, hence $g \in \ell.\mathrm{ann}_{R[x,\sigma]}(f_n)$. Therefore, $R[x,\sigma]$ is left endo-Noetherian.

\vspace{0.3cm}
\noindent
$(b)\Longrightarrow (a)$. Let $(q_i)_{i \in \mathbb{N}}$ be a sequence of elements of $T + R[x,\sigma]x$ such that:
\[
\ell.\mathrm{ann}_{T + R[x,\sigma]x}(q_1) \subseteq \ell.\mathrm{ann}_{T + R[x,\sigma]x}(q_2) \subseteq \cdots.
\]

We will show that:
\[
\ell.\mathrm{ann}_{R[x,\sigma]}(q_i) \subseteq \ell.\mathrm{ann}_{R[x,\sigma]}(q_{i+1}) \quad \text{for each } i \in \mathbb{N}.
\]

Let $g \in \ell.\mathrm{ann}_{R[x,\sigma]}(q_i)$. Then:
\[
g q_i = 0 \Rightarrow x g q_i = \sigma(g) x q_i = 0,
\]
so $\sigma(g) x \in \ell.\mathrm{ann}_{T + R[x,\sigma]x}(q_i) \subseteq \ell.\mathrm{ann}_{T + R[x,\sigma]x}(q_{i+1})$.

Hence:
\[
\sigma(g) x q_{i+1} = \sigma(g) \sigma(q_{i+1}) x = 0 \Rightarrow \sigma(g) \sigma(q_{i+1}) = 0.
\]
Again, since $\sigma$ is injective, we get $g q_{i+1} = 0$, i.e., $g \in \ell.\mathrm{ann}_{R[x,\sigma]}(q_{i+1})$.

Now, since $R[x,\sigma]$ is left endo-Noetherian, there exists a positive integer $n$ such that for all $k \geq n$:
\[
\ell.\mathrm{ann}_{R[x,\sigma]}(q_k) = \ell.\mathrm{ann}_{R[x,\sigma]}(q_n).
\]

We will now show that:
\[
\ell.\mathrm{ann}_{T + R[x,\sigma]x}(q_k) = \ell.\mathrm{ann}_{T + R[x,\sigma]x}(q_n)
\quad \text{for each } k \geq n.
\]

Let $q \in \ell.\mathrm{ann}_{T + R[x,\sigma]x}(q_k)$. Then:
\[
q q_k = 0 \Rightarrow q \in \ell.\mathrm{ann}_{R[x,\sigma]}(q_k) \subseteq \ell.\mathrm{ann}_{R[x,\sigma]}(q_n),
\]
so $q q_n = 0$, i.e., $q \in \ell.\mathrm{ann}_{T + R[x,\sigma]x}(q_n)$. Thus, $T + R[x,\sigma]x$ is left endo-Noetherian.

	\end{proof}
Similarly, we can deduce the following proposition: 
 \begin{proposition}
    Let $T \subseteq R$ be an extension of rings and $\sigma$ an injective endomorphism of $R$. Then the following conditions are equivalent:
    \begin{enumerate}
        \item $T+xR[x,\sigma]$ is right endo-Noetherian.
        \item $R[x,\sigma]$ is right endo-Noetherian.
    \end{enumerate}
   \end{proposition}
	
 According to \cite{hong2003skew}, a ring $R$ is called $\sigma$-skew Armendariz if $f(x)g(x) = 0$ for $f(x)=\sum_{i=0}^{n} a_{i} x^{i} $and $g(x)=\sum_{j=0}^{m} b_{j}  x^{j} \in  R[x,\sigma] $, then $a_{i}\sigma^{i}(b_{j} ) = 0$ for all $i$, $j$. 

On the other hand, we assume that $R$ is $\sigma$-skew Armendariz and $\sigma$-compatible in order for the structures $T+R[x,\sigma]x$ and $T+xR[x,\sigma]$ to be right endo-Noetherian and left endo-Noetherian, respectively.
\begin{lemma} 
     Let $R$ be a ring, $ \sigma$ an endomorphism of $R$ and $R$ $\sigma$-skew Armendariz and $\sigma$-compatible. Then for every two polynomials $f(x)=\sum_{i=0}^{n} a_{i} x^{i} $and $g(x)=\sum_{j=0}^{m} b_{j}  x^{j}$, $f(x) g(x)=0$ in $ R[x,\sigma] $ if and only if  $f(x) \sigma (g(x))=0$. 
\end{lemma}
\begin{proof}
    Let $f(x)=\sum_{i=0}^{n} a_{i} x^{i} $and $g(x)=\sum_{j=0}^{m} b_{j}  x^{j}$
    
$\Longleftrightarrow f(x) g(x)=\sum_{k=0}^{n+m}c_{k}  x^{k}=0$, $c_{k}=\sum_{i+j=k} a_{i} \sigma^{i} (b_j)$
  
$\Longleftrightarrow \sum_{i+j=k}a_{i} \sigma^{i} (b_{j})=0$, since $ R$ 
is $\sigma$-skew Armendariz, we have $a_i \sigma^{i} (b_{j})=0$ for all $0\leq i\leq n$, $0\leq j\leq m$, 
and since $R$ is $\sigma$-compatible, we have $a_{i} \sigma(\sigma^{i} (b_{j}))=a_{i} \sigma^{i+1} (b_{j})=0 $ for all $0\leq i\leq n$, $0\leq j\leq m 
\Longleftrightarrow c^{'}_k=\sum_{i+j=k}a_{i} \sigma^{i+1} (b_{j})=0
\Longleftrightarrow f(x) \sigma(g(x))=\sum_{k=0}^{n+m}c^{'}_{k}  x^{k} =0$.

\end{proof}

	\begin{theorem}
	    Let $T \subseteq R$ be an extension of rings and $\sigma$ be an endomorphism of $R$. Assume that $R$ is $\sigma$-skew Armendariz and $\sigma$-compatible. Then the following conditions are equivalent:
     \begin{enumerate}
         \item $T+R[x,\sigma]x$ is right endo-Noetherian.
         \item $R[x,\sigma]$ is right endo-Noetherian.
     \end{enumerate}
	
	\end{theorem}
\begin{proof}
\noindent
$(a)\Longrightarrow (b)$. Let $(f_k)_{k \in \mathbb{N}}$ be a sequence of elements in $R[x,\sigma]$ such that
\[
r.\mathrm{ann}_{R[x,\sigma]}(f_1) \subseteq r.\mathrm{ann}_{R[x,\sigma]}(f_2) \subseteq \cdots.
\]
Since $f_k x \in T + R[x,\sigma]x$ for each $f_k \in R[x,\sigma]$, we will show that
\[
r.\mathrm{ann}_{T + R[x,\sigma]x}(f_i x) \subseteq r.\mathrm{ann}_{T + R[x,\sigma]x}(f_{i+1} x) \quad \text{for each } i \in \mathbb{N}.
\]

Let $q \in r.\mathrm{ann}_{T + R[x,\sigma]x}(f_i x)$. Then:
\[
f_i \sigma(q) = f_i \sigma(q)x = f_i x q = 0,
\]
so $\sigma(q) \in r.\mathrm{ann}_{R[x,\sigma]}(f_i) \subseteq r.\mathrm{ann}_{R[x,\sigma]}(f_{i+1})$. Hence:
\[
f_{i+1} x q = f_{i+1} \sigma(q)x = f_{i+1} \sigma(q) = 0,
\]
and thus $q \in r.\mathrm{ann}_{T + R[x,\sigma]x}(f_{i+1} x)$.

Now, since $T + R[x,\sigma]x$ is right endo-Noetherian, there exists a positive integer $n$ such that for each $k \geq n$:
\[
r.\mathrm{ann}_{T + R[x,\sigma]x}(f_k x) = r.\mathrm{ann}_{T + R[x,\sigma]x}(f_n x).
\]

We will now show that:
\[
r.\mathrm{ann}_{R[x,\sigma]}(f_k) = r.\mathrm{ann}_{R[x,\sigma]}(f_n) \quad \text{for each } k \geq n.
\]

Let $g \in r.\mathrm{ann}_{R[x,\sigma]}(f_k)$. Then $f_k g = 0$. Since $R$ is $\sigma$-skew Armendariz and $\sigma$-compatible (by Lemma 2), we get:
\[
f_k \sigma(g) = 0 \quad \Rightarrow \quad f_k x g x = f_k \sigma(g) x^2 = 0,
\]
so $g x \in r.\mathrm{ann}_{T + R[x,\sigma]x}(f_k x)$.

Therefore:
\[
f_n \sigma(g) = f_n \sigma(g) x^2 = f_n x g x = 0.
\]
By the same lemma, this implies $f_n g = 0$, so $g \in r.\mathrm{ann}_{R[x,\sigma]}(f_n)$. Hence, $R[x,\sigma]$ is right endo-Noetherian.

\vspace{0.3cm}
\noindent
$(b)\Longrightarrow (a)$. Note that this implication always holds and does not require the assumption that $R$ is $\sigma$-skew Armendariz or $\sigma$-compatible.

Let $(q_i)_{i \in \mathbb{N}}$ be a sequence of elements in $T + R[x,\sigma]x$ such that:
\[
r.\mathrm{ann}_{T + R[x,\sigma]x}(q_1) \subseteq r.\mathrm{ann}_{T + R[x,\sigma]x}(q_2) \subseteq \cdots.
\]

We will show that:
\[
r.\mathrm{ann}_{R[x,\sigma]}(q_i) \subseteq r.\mathrm{ann}_{R[x,\sigma]}(q_{i+1}) \quad \text{for each } i \in \mathbb{N}.
\]

Let $g \in r.\mathrm{ann}_{R[x,\sigma]}(q_i)$. Then:
\[
q_i g x = q_i g = 0 \quad \Rightarrow \quad g x \in r.\mathrm{ann}_{T + R[x,\sigma]x}(q_i),
\]
and so:
\[
q_{i+1} g = q_{i+1} g x = 0 \quad \Rightarrow \quad g \in r.\mathrm{ann}_{R[x,\sigma]}(q_{i+1}).
\]

Since $R[x,\sigma]$ is right endo-Noetherian, there exists a positive integer $n$ such that for all $k \geq n$:
\[
r.\mathrm{ann}_{R[x,\sigma]}(q_k) = r.\mathrm{ann}_{R[x,\sigma]}(q_n).
\]

We will now show that:
\[
r.\mathrm{ann}_{T + R[x,\sigma]x}(q_k) = r.\mathrm{ann}_{T + R[x,\sigma]x}(q_n) \quad \text{for each } k \geq n.
\]

Let $q \in r.\mathrm{ann}_{T + R[x,\sigma]x}(q_k)$. Then $q_k q = 0$ and $q \in r.\mathrm{ann}_{R[x,\sigma]}(q_k) \subseteq r.\mathrm{ann}_{R[x,\sigma]}(q_n)$. Hence:
\[
q_n q = 0 \quad \Rightarrow \quad q \in r.\mathrm{ann}_{T + R[x,\sigma]x}(q_n).
\]

Therefore, $T + R[x,\sigma]x$ is right endo-Noetherian.

\end{proof}
Similarly, we can deduce the following proposition:
\begin{proposition}
    Let $T \subseteq R$ be an extension of rings and $\sigma$ an endomorphism of $R$. Assume that $R$ is $\sigma$-skew Armendariz and $\sigma$-compatible. Then the following conditions are equivalent:
    \begin{enumerate}
        \item $T+xR[x,\sigma]$ is left endo-Noetherian.
        \item $R[x,\sigma]$ is left endo-Noetherian.
    \end{enumerate}
\end{proposition}

\section{More Results on Endo-Noetherian Rings}
Let $S=(R_{n})_{n\in \mathbb{N}}$ be an increasing sequence of rings, $ R=\cup _{n\in \mathbb{N}} R_{n} $, and $S[x]$ the ring of polynomials with coefficients of degree $n$ in $R_{n}$. In \cite[Theorem 2.1]
{hamed2021rings} the authors proved in commutative case that the ring $R$ is strongly Hopfian if
and only if its polynomial $ R[x]$ is strongly Hopfian. In the following, we generalize this theorem to the noncommutative case.
\begin{theorem}
  Let $S=(R_{n})_{n\in \mathbb{N}}$ be an increasing sequence of rings, $ R=\cup _{n\in \mathbb{N}} R_{n} $, and $S[x]$ the ring of polynomials with coefficients of degree $n$ in $R_{n}$. The following conditions are equivalent:
 \begin{enumerate}
     \item $S[x]$ is left endo-Noetherian.
     \item $R[x]$ is left endo-Noetherian.
 \end{enumerate}
 \end{theorem}
	\begin{proof}
\noindent
$(a) \Longrightarrow (b)$. Let $(f_i(x))_{i \in \mathbb{N}}$ be a sequence of elements of $R[x]$ such that
\[
\ell.\mathrm{ann}_{R[x]}(f_1(x)) \subseteq \ell.\mathrm{ann}_{R[x]}(f_2(x)) \subseteq \cdots.
\]

Note that if $f_i(x) = \sum_{j_i=0}^{m_i} a_{j_i} x^{j_i} \in R[x]$, then for each $j_i$, $0 \leq j_i \leq m_i$, there exists $t_{j_i} \in \mathbb{N}$ such that $a_{j_i} \in R_{t_{j_i}}$. Let
\[
l_i = \max \{ t_{j_i} \mid 0 \leq j_i \leq m_i \}.
\]
Then $f_i(x)x^{l_i} \in S[x]$ for each $i \in \mathbb{N}$. We will show that
\[
\ell.\mathrm{ann}_{S[x]}(f_i(x) x^{l_i}) \subseteq \ell.\mathrm{ann}_{S[x]}(f_{i+1}(x) x^{l_{i+1}}).
\]

Let $g(x) \in \ell.\mathrm{ann}_{S[x]}(f_i(x) x^{l_i})$. Then
\[
g(x)f_i(x) = g(x)f_i(x)x^{l_i} = 0,
\]
so $g(x) \in \ell.\mathrm{ann}_{R[x]}(f_i(x)) \subseteq \ell.\mathrm{ann}_{R[x]}(f_{i+1}(x))$.

Thus,
\[
g(x)f_{i+1}(x) = 0 \quad \text{and} \quad g(x)f_{i+1}(x)x^{l_{i+1}} = 0,
\]
hence $g(x) \in \ell.\mathrm{ann}_{S[x]}(f_{i+1}(x) x^{l_{i+1}})$.

Now, since $S[x]$ is left endo-Noetherian, there exists a positive integer $n$ such that for all $k \geq n$,
\[
\ell.\mathrm{ann}_{S[x]}(f_k(x)x^{l_k}) = \ell.\mathrm{ann}_{S[x]}(f_n(x)x^{l_n}).
\]

We will show that for each $k \geq n$,
\[
\ell.\mathrm{ann}_{R[x]}(f_k(x)) = \ell.\mathrm{ann}_{R[x]}(f_n(x)).
\]

Let $h(x) \in \ell.\mathrm{ann}_{R[x]}(f_k(x))$, so $h(x)f_k(x) = 0$. Then $h(x)f_k(x)x^{l_k} = 0$. As above, there exists $m \in \mathbb{N}$ such that $h(x)x^m \in S[x]$. Therefore,
\[
h(x)x^m f_k(x) x^{l_k} = 0 \quad \Rightarrow \quad h(x)x^m \in \ell.\mathrm{ann}_{S[x]}(f_k(x)x^{l_k}).
\]

So,
\[
h(x)x^m f_n(x)x^{l_n} = 0 \quad \Rightarrow \quad h(x)f_n(x) = 0,
\]
which means $h(x) \in \ell.\mathrm{ann}_{R[x]}(f_n(x))$.

Hence, $R[x]$ is left endo-Noetherian.

\vspace{0.3cm}
\noindent
$(b) \Longrightarrow (a)$. Let $(q_i(x))_{i \in \mathbb{N}}$ be a sequence of elements in $S[x]$ such that
\[
\ell.\mathrm{ann}_{S[x]}(q_1(x)) \subseteq \ell.\mathrm{ann}_{S[x]}(q_2(x)) \subseteq \cdots.
\]

We will show that
\[
\ell.\mathrm{ann}_{R[x]}(q_i(x)) \subseteq \ell.\mathrm{ann}_{R[x]}(q_{i+1}(x)) \quad \text{for each } i \in \mathbb{N}.
\]

Let $g(x) \in \ell.\mathrm{ann}_{R[x]}(q_i(x))$. Then $g(x)q_i(x) = 0$. As above, there exists $s \in \mathbb{N}$ such that $x^s g(x) \in S[x]$ and
\[
x^s g(x) q_i(x) = 0 \quad \Rightarrow \quad x^s g(x) \in \ell.\mathrm{ann}_{S[x]}(q_i(x)) \subseteq \ell.\mathrm{ann}_{S[x]}(q_{i+1}(x)).
\]

Hence,
\[
g(x)q_{i+1}(x) = x^s g(x) q_{i+1}(x) = 0 \quad \Rightarrow \quad g(x) \in \ell.\mathrm{ann}_{R[x]}(q_{i+1}(x)).
\]

Since $R[x]$ is left endo-Noetherian, there exists a positive integer $n$ such that for all $k \geq n$,
\[
\ell.\mathrm{ann}_{R[x]}(q_k(x)) = \ell.\mathrm{ann}_{R[x]}(q_n(x)).
\]

We will now show that
\[
\ell.\mathrm{ann}_{S[x]}(q_k(x)) = \ell.\mathrm{ann}_{S[x]}(q_n(x)) \quad \text{for all } k \geq n.
\]

Let $h(x) \in \ell.\mathrm{ann}_{S[x]}(q_k(x))$. Then $h(x)q_k(x) = 0$, and since $h(x) \in \ell.\mathrm{ann}_{R[x]}(q_k(x)) = \ell.\mathrm{ann}_{R[x]}(q_n(x))$, we get:
\[
h(x)q_n(x) = 0 \quad \Rightarrow \quad h(x) \in \ell.\mathrm{ann}_{S[x]}(q_n(x)).
\]

Therefore, $S[x]$ is left endo-Noetherian.

	\end{proof}

\begin{proposition}\cite[Corollary 2.1]{mohamed2023endo} \label{prop2}.
    Let $R$ be an Armendariz ring. Then, the following statements are equivalent:
        \begin{enumerate}
    	\item $R[x]$ is left endo-Noetherian.
    	\item $R$ satisfies the acc on left annihilators of finite subset.
    	\item $R$ satisfies the acc on left annihilators of finitely generated ideals of $R$.
    	\end{enumerate}
    	In particular, if $R[x]$ is left endo-Noetherian, then $R$ is left endo-Noetherian.
\end{proposition}
\begin{corollary}
     Let $S$, $R$ and $S[x]$  be as in Theorem 6. 
If $S[x]$ is left endo-Noetherian, then $R$ is left endo-Noetherian.
\end{corollary}
\begin{proof}
Assume that $S[x]$ is left endo-Noetherian. From Theorem 6, $R[x]$ is left endo-Noetherian and from Proposition 7, $R$ is left endo-Noetherian.
\end{proof}

\begin{corollary}
    Let $S$, $R$ and $S[x]$  be as in Theorem 6. If $R$ satisfies the acc on left annihilators of finite subset, then $S[x]$ is left endo-Noetherian.
\end{corollary}
\begin{proof}
Let $R$ be an Armendariz ring satisfies the acc on left annihilators of finite subset. From Proposition 7, $R[x]$ is left endo-Noetherian, and from Theorem 6, $S[x]$ is left endo-Noetherian. 
\end{proof}
In the following, we study when the amalgamated rings are left endo-Noetherian. 

In \cite{d2007amalgamated}, M. D’Anna and M. Fontana introduced a construction called the amalgamated
duplication of a ring $R$ along an ideal $I$ of $R$, denoted by $ R \bowtie I $, it is defined as the following subring of $R \times R$:\begin{equation*}
 R \bowtie I = \{(r , r + i ) \mid r \in R, i \in I \}.
\end{equation*}

Recall that, in \cite{gouaid2020endo} an ideal $I$ of a ring $R$ is called a regular ideal if it contains a regular element.
In the following theorem, we give the necessary and
sufficient conditions for the ring $ R \bowtie I $ to be left endo-Noetherian.

\begin{theorem}
     Let $R$ be a ring and $I$ a regular ideal of $R$. Then, the following
assertions are equivalent:
\begin{enumerate}
    \item $R$ is left endo-Noetherian.
    \item $R \times R$ is left endo-Noetherian.
    \item $R \bowtie I$ is left endo-Noetherian.
\end{enumerate}

\end{theorem}
\begin{proof}
\noindent
$(a)\Longrightarrow (b)$. It follows from \cite[Theorem 2]{mohamed2023endo}.

\vspace{0.3cm}
\noindent
$(b)\Longrightarrow (c)$. Let 
\[
(r_1, r_1+i_1),\ (r_2, r_2+i_2),\ \dots \in R \bowtie I
\]
such that
\[
\ell.\mathrm{ann}_{R \bowtie I}(r_1, r_1+i_1) \subseteq \ell.\mathrm{ann}_{R \bowtie I}(r_2, r_2+i_2) \subseteq \cdots.
\]

We will show that
\[
\ell.\mathrm{ann}_{R \times R}(r_k, r_k+i_k) \subseteq \ell.\mathrm{ann}_{R \times R}(r_{k+1}, r_{k+1}+i_{k+1}) \quad \text{for each } k \geq 1.
\]

Let $(a, b) \in \ell.\mathrm{ann}_{R \times R}(r_k, r_k+i_k)$, and let $i$ be a regular element of $I$. Consider the element $(i, i) \in R \bowtie I$. Then,
\[
(i,i)(a,b) = (ia, ia + i(b-a)) \in R \bowtie I.
\]

Moreover,
\[
(i,i)(a,b)(r_k, r_k+i_k) = (0, 0),
\]
so $(i,i)(a,b) \in \ell.\mathrm{ann}_{R \bowtie I}(r_{k+1}, r_{k+1}+i_{k+1})$. Since $(i,i)$ is regular in $R \times R$, we conclude that
\[
(a, b) \in \ell.\mathrm{ann}_{R \times R}(r_{k+1}, r_{k+1}+i_{k+1}).
\]

By the hypothesis that $R \times R$ is left endo-Noetherian, there exists a positive integer $n$ such that for all $k \geq n$,
\[
\ell.\mathrm{ann}_{R \times R}(r_k, r_k+i_k) = \ell.\mathrm{ann}_{R \times R}(r_n, r_n+i_n).
\]

Thus, for each $k \geq n$,
\[
(R \bowtie I) \cap \ell.\mathrm{ann}_{R \times R}(r_k, r_k+i_k) = (R \bowtie I) \cap \ell.\mathrm{ann}_{R \times R}(r_n, r_n+i_n),
\]
which implies
\[
\ell.\mathrm{ann}_{R \bowtie I}(r_k, r_k+i_k) = \ell.\mathrm{ann}_{R \bowtie I}(r_n, r_n+i_n).
\]

\vspace{0.3cm}
\noindent
$(c)\Longrightarrow (a)$. Note that this implication is always true and does not require the assumption that $I$ contains a regular element.

Let $r_1, r_2, \dots \in R$ such that
\[
\ell.\mathrm{ann}_R(r_1) \subseteq \ell.\mathrm{ann}_R(r_2) \subseteq \cdots.
\]

We will show that
\[
\ell.\mathrm{ann}_{R \bowtie I}(r_k, r_k) \subseteq \ell.\mathrm{ann}_{R \bowtie I}(r_{k+1}, r_{k+1}) \quad \text{for each } k \geq 1.
\]

Let $(\alpha, \alpha+i) \in \ell.\mathrm{ann}_{R \bowtie I}(r_k, r_k)$. Then $\alpha, \alpha+i \in \ell.\mathrm{ann}_R(r_k)$, so
\[
(\alpha, \alpha+i) \in \ell.\mathrm{ann}_{R \bowtie I}(r_{k+1}, r_{k+1}).
\]

Since $R \bowtie I$ is left endo-Noetherian, there exists a positive integer $n$ such that for all $k \geq n$,
\[
\ell.\mathrm{ann}_{R \bowtie I}(r_k, r_k) = \ell.\mathrm{ann}_{R \bowtie I}(r_n, r_n).
\]

We now show that
\[
\ell.\mathrm{ann}_R(r_k) = \ell.\mathrm{ann}_R(r_n) \quad \text{for all } k \geq n.
\]

Let $b \in \ell.\mathrm{ann}_R(r_k)$. Then, since $(b, b) \in \ell.\mathrm{ann}_{R \bowtie I}(r_k, r_k)$ and $(b, b)(r_n, r_n) = (0, 0)$, we conclude that $br_n = 0$. Therefore,
\[
b \in \ell.\mathrm{ann}_R(r_n).
\]

Hence, $R$ is left endo-Noetherian.

\end{proof}
In \cite{d2009amalgamated}, M. D’Anna and M. Fontana introduced a new ring construction of amalgamated algebra called the amalgamation of $ R$ with $S$ along $J$ with respect to $f$, denoted by $R \bowtie ^f J$, as a generalization of the amalgamated
duplication $ R \bowtie I $, it is defined as the following subring of $R \times S$:
\begin{equation*}
     R \bowtie ^f J:=\{(r, f (r)+ j) \mid r \in R, j \in J\}
\end{equation*}
for a given ring homomorphism $f : R \longrightarrow S$ and ideal $J$ of $S$.\\

In the next proposition we show when $ R \bowtie ^f J$ is left endo-Noetherian.
\begin{proposition}
     Let $R$ and $S$ be two rings, $J$ be an ideal of $S$, and let $f : R \longrightarrow S$  be a ring homomorphism. If $R$ and $f(R)+J$ are left endo-Noetherian, then $ R \bowtie ^f J$ is left endo-Noetherian.
\end{proposition}
\begin{proof}
    Let $(r_{i},f(r_{i} )+j_{i})_{i\in \mathbb{N}}$ be a sequence of elements of $ R \bowtie ^f J$ such that\\ $\ell .\mathrm {ann}_{ R \bowtie ^f J} (r_1,f(r_1 )+j_1 )\subseteq  \ell .\mathrm {ann}_{ R \bowtie ^f J} (r_2,f(r_2 )+j_2 )\subseteq \dots$ . Since $R$ and $f(R)+J$ are left endo-Noetherian, there exists a positive integer $n$ such that $\ell .\mathrm {ann}_R (r_k )=\ell .\mathrm {ann}_{R} (r_n )$ and $\ell .\mathrm {ann}_{f(R)+J}(f(r_k )+j_k )= \ell .\mathrm {ann}_{f(R)+J} (f(r_n )+j_n )$ for each $k\geq n$ . Hence $\ell .\mathrm {ann}_{ R \bowtie ^f J} (r_k,f(r_k )+j_k )=  \ell .\mathrm {ann}_{ R \bowtie ^f J} (r_n,f(r_n )+j_n )$, and $ R \bowtie ^f J$  is left endo-Noetherian.

\end{proof}

\bibliographystyle{splncs04}
\bibliography{biblio}

\begin{thebibliography}{10}
\providecommand{\url}[1]{\texttt{#1}}
\providecommand{\urlprefix}{URL }
\providecommand{\doi}[1]{https://doi.org/#1}

\bibitem{ahmed2015s}
Ahmed, H., Sana, H.: S-noetherian rings of the forms $\mathbb{A}[x]$ and $\mathbb{A}[[x]]$. Communications in Algebra  \textbf{43}(9),  3848--3856 (2015)

\bibitem{d2007amalgamated}
D'Anna, M., Fontana, M.: An amalgamated duplication of a ring along an ideal: the basic properties. Journal of Algebra and its Applications  \textbf{6}(03),  443--459 (2007)

\bibitem{dina2024almost}
Dina~Abdelhakim, R. M.~Salem, S.E.D.: Almost right (left) semiclean rings of skew generalized power series. Journal of Scientific Research in Science  (2024)

\bibitem{d2009amalgamated}
D’Anna, M., Finocchiaro, C.A., Fontana, M.: Amalgamated algebras along an ideal. Commutative algebra and its applications pp. 155--172 (2009)

\bibitem{gouaid2020endo}
Gouaid, B., Hamed, A., Benhissi, A.: Endo-noetherian rings. Annali di Matematica Pura ed Applicata (1923-)  \textbf{199}(2),  563--572 (2020)

\bibitem{hamed2021rings}
Hamed, A., Gouaid, B., Benhissi, A.: Rings satisfying the strongly hopfian and s-strongly hopfian properties. Math. Rep  \textbf{23}(4),  383--395 (2021)

\bibitem{hong2003skew}
Hong, C.Y., Kim, N.K., Kwak, T.K.: On skew armendariz rings  (2003)

\bibitem{kaidimodules}
Kaidi, A.: Modules with chain conditions on endoimages endokernels. Preprint

\bibitem{mcconnell2001noncommutative}
McConnell, J.C., Robson, J.C., Small, L.W.: Noncommutative noetherian rings, vol.~30. American Mathematical Soc. (2001)

\bibitem{mohamed2023endo}
Mohamed, N.A., Salem, R.M., Abdel-Khalek, R.E.: Endo-noetherian skew generalized power series rings. Assiut University Journal of Multidisciplinary Scientific Research  \textbf{52}(1),  13--22 (2023)

\end{thebibliography}

\end{document}